\documentclass[11pt]{article}
\usepackage{amsmath,amssymb,amsthm,latexsym,mathrsfs}
\usepackage{indentfirst}
 \DeclareMathOperator{\DIV}{div}
  \DeclareMathOperator{\grad}{grad}
 \DeclareMathOperator{\trace}{tr}

 \newcommand{\dif}{\mathrm{d}}
 \newcommand{\ROM}[1]{\mathrm{\uppercase\expandafter{\romannumeral#1}}}
  \theoremstyle{definition}
 \newtheorem{theorem}{Theorem}[section]
 \newtheorem{lemma}[theorem]{Lemma}

  \makeatletter
\renewcommand \thesection {\S\@arabic\c@section}
\renewcommand \thetheorem {\@arabic\c@section.\@arabic\c@theorem}
\makeatother
\title{\textbf{ Integral formula of Minkowski type and new characterization of the Wulff shape}}
\author {Yijun He
\thanks {Partially supported by Tianyuan Fund for Mathematics of NSFC.}
 \and Haizhong Li
\thanks {Partially supported by the grant No. 10531090 of the NSFC and by Doctoral Program Foundation of the Ministry of Education of China (2006).}}
\date{}
\begin{document}
\maketitle
\begin{abstract}
Given a positive function $F$ on $S^n$ which satisfies a convexity condition, we introduce the $r$-th anisotropic mean curvature $M_r$ for
hypersurfaces in $\mathbb{R}^{n+1}$ which is a generalization of the usual $r$-th mean curvature $H_r$. We get integral formulas of Minkowski
type for compact hypersurfaces in $R^{n+1}$. We give some new characterizations of the Wulff shape by use of our integral formulas of Minkowski
type, in case $F=1$ which reduces to some well-known results.
\end{abstract}

\section{Introduction}
Let $F\colon S^n\to\mathbb{R}^+$ be a smooth function which
satisfies the following convexity condition:
\begin{equation}
(D^2F+F1)_x>0,\quad\forall x\in S^n,
\end{equation}
where $D^2F$ denotes the intrinsic Hessian of $F$ on $S^n$ and 1
denotes the identity on $T_x S^n$, $>0$ means that the matrix is
positive definite. We consider the map
$$\phi\colon S^n\to\mathbb{R}^{n+1},$$
$$x\to F(x)x+(\grad_{S^n}F)_x,$$
its image $W_F=\phi(S^n)$ is a smooth, convex hypersurface in $\mathbb{R}^{n+1}$ called the Wulff shape of $F$ (see [3], [6], [10], [13], [14]).

Now let $X\colon M\to\mathbb{R}^{n+1}$ ba a smooth immersion of a compact, orientable hypersurface without boundary. Let $\nu\colon M\to S^n$
denotes its Gauss map, that is, $\nu$ is an unit inner normal vector of $M$.

Let $A_F=D^2F+F1$, $S_F=-A_F\circ\dif\nu$. $S_F$ is called the
$F$-Weingarten operator, and the eigenvalues of $S_F$ are called
anisotropic principal curvatures. Let $\sigma_r$ be the elementary
symmetric functions of the anisotropic principal curvatures
$\lambda_1, \lambda_2, \cdots, \lambda_n$:
$$\sigma_r=\sum_{i_1<\cdots<i_r}\lambda_{i_1}\cdots\lambda_{i_r}\quad (1\leq r\leq n).$$
We set $\sigma_0=1$. The $r$-anisotropic mean curvature $M_r$ is
defined by $M_r=\sigma_r/C^r_n$.

In this paper we first give the following integral formulas of Minkowski type for compact hypersurfaces in $\mathbb{R}^{n+1}$.
\begin{theorem}
Let $X\colon M\to\mathbb{R}^{n+1}$ be an $n$-dimensional compact hypersurface, $F\colon S^n\to\mathbb{R}^+$ be a smooth function which satisfies
(1), then we have the following integral formulas of Minkowski type hold:
\begin{equation}\label{Min}
 \int_M(FM_r+M_{r+1}\langle X, \nu\rangle)\dif A_X=0,\quad r=0, 1,
 \cdots, n-1.
\end{equation}
\end{theorem}
By use of above integral formulas of Minkowski type, we prove the following new characterizations of the Wulff shape:
\begin{theorem} Let $X\colon M\to\mathbb{R}^{n+1}$ be an $n$-dimensional compact hypersurface, $F\colon S^n\to\mathbb{R}^+$ be a smooth function
which satisfies (1), and $M_1=\mbox{const}$ and $\langle X, \nu\rangle$ has fixed sign, then up to translations and homotheties, $X(M)$ is the
Wulff  shape.
\end{theorem}
\begin{theorem}Let $X\colon M\to\mathbb{R}^{n+1}$ be an $n$-dimensional compact hypersurface, $F\colon S^n\to\mathbb{R}^+$ be a smooth
function which satisfies (1).
 If $M_1=\mbox{const}$ and $M_r=\mbox{const}$ for some $r$, $2\leq r\leq n$, then up to translations and homotheties, $X(M)$ is the Wulff shape.
\end{theorem}
\begin{theorem}Let $X\colon M\to\mathbb{R}^{n+1}$ be an $n$-dimensional compact convex hypersurface, $F\colon S^n\to\mathbb{R}^+$ be a smooth
function which satisfies (1).
 If $\frac{M_r}{M_k}=\mbox{const}$ for some $k$ and $r$, with $0\leq k<r\leq n$, then then up to translations and homotheties,
 $X(M)$ is the Wulff shape.
\end{theorem}
\begin{theorem}Let $X\colon M\to\mathbb{R}^{n+1}$ be an $n$-dimensional compact hypersurface, $F\colon S^n\to\mathbb{R}^+$ be a smooth
function which satisfies (1).
 If $\frac{M_k}{M_n}=\mbox{const}$ for some $k$, with $0\leq k\leq n-1$, then then up to translations and homotheties,
 $X(M)$ is the Wulff shape.
\end{theorem}
Choosing $k=0$ in Theorem 1.4, we get

\medskip\noindent
{\bf Corollary 1.1} {\it Let $X\colon M\to\mathbb{R}^{n+1}$ be an $n$-dimensional compact convex hypersurface, $F\colon S^n\to\mathbb{R}^+$ be a
smooth function which satisfies (1),
 and for a fixed $r$ with $1\leq r\leq n$, $M_r=\mbox{const}$, then up to translations and homotheties, $X(M)$ is the Wulff shape.}

\medskip\noindent
{\bf Remark 1.1}  When $F=1$, Wulff shape is just the round sphere and $M_r=H_r$, formula (2) reduces to the classical Minkowski integral
  formula (see [5] or [11]), Theorem 1.2 reduces to the classical Theorem given by S\"{u}ss [12], Corollary 1.1 reduces to Theorem of
  Yano [15],  Theorem 1.3 reduces to Theorem of Choe [2].

\section{Preliminaries}

Let $\{E_1,\cdots,E_n\}$ is a local orthogonal frame on $S^n$, let $e_i=E_i\circ\nu$, where $i=1,\cdots,n$, then $\{e_1,\cdots,e_n\}$ is a local
orthogonal frame of $X\colon M\to\mathbb{R}^{n+1}$.

The structure equation of $S^n$ is:
$$\left\{
\begin{array}
  {l}
  \dif x=\sum_i\theta_iE_i\\
  \dif E_i=\sum_j\theta_{ij}E_j-\theta_ix\\
  \dif\theta_i=\sum_j\theta_{ij}\wedge\theta_j\\
  \dif\theta_{ij}-\sum_k\theta_{ik}\wedge\theta_{kj}=-\frac12\sum_{kl}{\tilde R}_{ijkl}\theta_k\wedge\theta_l=-\theta_i\wedge\theta_j
\end{array}\right.
\eqno (3)
$$
where $\theta_{ij}+\theta_{ji}=0$ and
$$
 {\tilde R}_{ijkl}=\delta_{ik}\delta_{jl}-\delta_{il}\delta_{jk}.
\eqno (4)
$$

The structure equation of $X$ is (see [7], [8]):
$$\left\{
\begin{array}
  {l}
  \dif X=\sum_i\omega_ie_i\\
  \dif \nu=-\sum_{ij}h_{ij}\omega_je_i\\
  \dif
  e_i=\sum_j\omega_{ij}e_j+\sum_jh_{ij}\omega_j\nu\\
  \dif\omega_i=\sum_j\omega_{ij}\wedge\omega_j\\
  \dif\omega_{ij}-\sum_k\omega_{ik}\wedge\omega_{kj}=
  -\frac12\sum_{kl}R_{ijkl}\theta_k\wedge\theta_l
\end{array}\right.
\eqno (5)
$$
where $\omega_{ij}+\omega_{ji}=0$, $R_{ijkl}+R_{ijlk}=0$, and $R_{ijkl}$ are the components of the Riemannian curvature tensor of $M$ with
respect to the induced metric $dX\cdot dX$.

From $\dif e_i=\dif(E_i\circ\nu)=\nu^*\dif
E_i=\sum_j\nu^*\theta_{ij}e_j-\nu^*\theta_i\nu$, we get
$$\left\{
\begin{array}
  {l}
  \omega_{ij}=\nu^*\theta_{ij}\\
  \nu^*\theta_i=-\sum_jh_{ij}\omega_j,
\end{array}\right.
\eqno (6)
$$
where $\omega_{ij}+\omega_{ji}=0$, $h_{ij}=h_{ji}$.

 Let $F\colon S^n\to\mathbb{R}^+$ be a smooth function, we denote the coefficients of covariant differential of $F$, $\grad_{S^n}F$,
$D^2F$ with respect to $\{E_i\}_{i=1,\cdots,n}$ by $F_i, F_{ij}, F_{ijk}$ respectively.

From Ricci identity and (4), we have
$$
F_{ijk}-F_{ikj}=\sum_m
 {\tilde R}_{mijk}=\delta_{ik}F_j-\delta_{ij}F_k,
\eqno (7)
$$
where $F_{ijk}$ denote the coefficients of the covariant differential of $F_{ij}$ on $S^n$.

 So, if we denote the coefficients of $A_F$ by $A_{ij}$, then we have from (7)
$$A_{ijk}=A_{jik}=A_{ikj},
\eqno (8)
$$
where $A_{ijk}$ denotes coefficient of the covariant differential of
$A_F$ on $S^n$.

Let $s_{ij}=\sum_k(A_{ik}\circ\nu) h_{kj}$, $S_F=-A_F\circ\dif\nu$, then we have $S_F(e_j)=\sum_is_{ij}e_i$. We call $S_F$ to be $F$-Weingarten
operator. From the positive definite of $(A_{ij})$ and the symmetry of $(h_{ij})$, we know the eigenvalues of $(s_{ij})$ are all real (in fact,
because $A=(A_{ij})$ is positive definite, there exists a non-singular matrix $C$ such that $A=C^tC$, we have $S=(s_{ij})=AB$ has the same
eigenvalues with the real symmetric matrix $CBC^T$, which follows from $|\lambda I-S|=|\lambda I-AB|=|\lambda I-C^tCB|=|\lambda I-CBC^t|$, where
$B=(h_{ij})$). we call them anisotropic principal curvatures, and denote them by $\lambda_1, \cdots, \lambda_n$.

We have $n$ invariants, the elementary symmetric function $\sigma_r$
of the anisotropic principal curvatures:
$$
  \sigma_r=\sum_{i_1<\cdots i_r}\lambda_{i_1}\cdots\lambda_{i_n}
  \quad (1\leq r\leq n).
\eqno (9)
$$
For convenience, we set $\sigma_0=1$. The $r$-anisotropic mean curvature $M_r$ is defined by
$$M_r=\sigma_r/C_n^r,\quad
C^r_n=\frac{n!}{r!(n-r)!}.
\eqno (10)
$$

Using the characteristic polynomial of $S_F$, $\sigma_r$ is defined
by
 $$
 \det(tI-S_F)=\sum_{r=0}^n(-1)^r\sigma_rt^{n-r}.
\eqno (11)
 $$
 So, we have
 $$
 \label{sigma}
   \sigma_r=\frac{1}{r!}\sum_{i_1,\cdots,i_r;j_1,\cdots,j_r}\delta_{i_1\cdots
   i_r}^{j_1\cdots j_r}s_{i_1j_1}\cdots s_{i_rj_r},
 \eqno (12)
 $$
  where $\delta_{i_1\cdots
   i_r}^{j_1\cdots j_r}$ is the usual generalized Kronecker symbol,
   i.e., $\delta_{i_1\cdots
   i_r}^{j_1\cdots j_r}$ equals +1 (resp. -1) if $i_1\cdots
   i_r$ are distinct and  $(j_1\cdots j_r)$
   is an even (resp. odd) permutation of $(i_1\cdots i_r)$ and in
   other cases it equals zero.

 We define $(F\circ\nu)_i, (F_i\circ\nu)_j, (A_{ij}\circ\nu)_k$ by
$$
\dif(F\circ\nu)=\sum_i(F\circ\nu)_i\omega_i,
\eqno (13)
$$
$$
\dif(F_i\circ\nu)+\sum_j(F_j\circ\nu)\omega_{ji}=\sum_j(F_i\circ\nu)_j\omega_j,
\eqno (14)
$$
$$
\dif(A_{ij}\circ\nu)+\sum(A_{kj}\circ\nu)\omega_{ki}+\sum_k(A_{ik}\circ\nu)\omega_{kj}= \sum_k(A_{ij}\circ\nu)_k\omega_k.
\eqno (15)
$$

By use of (3), (5) and (6), we have by a direct calculation
$$ \label{jb} \left\{\begin{array}
  {c}
(F\circ\nu)_i=-\sum_jh_{ij}F_j\circ\nu,
\\
(F_i\circ\nu)_j=-\sum_kh_{jk}F_{ik}\circ\nu,
\\
(A_{ij}\circ\nu)_k=-\sum_lh_{kl}A_{ijl}\circ\nu.
\end{array}\right.
\eqno (16)
$$

 \section{Some Lemmas}
  We introduce an important operator $P_r$ by
   $$
P_r=\sigma_rI-\sigma_{r-1}S_F+\cdots+(-1)^rS_F^r, \quad r=0, 1, \cdots, n.
 \eqno (17)
$$

  We have the following lemmas:
   \begin{lemma}\label{lemma3.1}
     $(S_FA_F)^t=S_FA_F$, $(\dif\nu\circ S_F)^t=\dif\nu\circ S_F$, $s_{ijk}=s_{ikj}$,
     $\sum_lh_{il}s_{lk}=\sum_lh_{kl}s_{li}$,
      $\sum_lh_{kl}(P_r)_{lj}=\sum_lh_{jl}(P_r)_{lk}$,
 where  $s_{ijk}$ are the components of the covariant derivarive of $s_{ij}$.
   \end{lemma}
   \begin{proof}
     Since $S_F=-A_F\circ\dif\nu$, and $A_F$, $\dif\nu$ are
     symmetric operators, the first two identities are obvious. From the symmetry property (8) of $A_{ijk}$, $h_{ij}=h_{ji}$ and Codazzi
     equation $h_{ijk}=h_{ikj}$, we have by use of (16)
     $$\begin{array}
       {rcl}
       s_{ijk} &=& (\sum_lA_{il} h_{lj})_k=\sum_l(A_{il}\circ\nu)_kh_{lj}+\sum_lA_{il}h_{ljk}\\
       &=&
     -\sum_{l,m}(A_{ilm}\circ\nu) h_{lj}h_{km}+\sum_{l}A_{il}h_{ljk}\\
     &=&
     \sum_m(A_{im}\circ\nu)_jh_{mk}+\sum_{l}A_{il} h_{lkj}=(\sum_lA_{il}h_{lk})_j=s_{ikj}.
   \end{array}
   \eqno (18)
  $$
$$
\sum_lh_{il}s_{lk}=\sum_{l,m}h_{il}A_{lm}h_{mk}=
 \sum_{l,m}h_{km}A_{ml}h_{li}=
\sum_lh_{kl}s_{li}.
$$

By use of the above formula and the definition of $P_r$, we get the last identity in Lemma 3.1.
 \end{proof}

   \begin{lemma}\label{lemma3.2}
   The matrix of $P_r$ is given by:
 $$
 \label{Pr}
 (P_r)_{ij}=\displaystyle{\frac{1}{r!}\sum_{i_1,\cdots,i_r;j_1,\cdots,j_r}
     \delta_{i_1\cdots i_r j}^{j_1\cdots j_r i}}s_{i_1j_1}\cdots s_{i_rj_r}
 \eqno (19)
 $$
 \end{lemma}
   \begin{proof}
     We prove Lemma 3.2 inductively. For $r=0$, it is easy to check that (19) is true.

     We can check directly
 $$
 \label{9}
 \delta_{i_1\cdots i_q}^{j_1\cdots j_q}=\left|
 \begin{array}{ccccc}
 \delta_{i_1}^{j_1} & \delta_{i_1}^{j_2} & \cdots & \delta_{i_1}^{j_{q-1}} & \delta_{i_1}^{j_q} \\
 \delta_{i_2}^{j_1} & \delta_{i_2}^{j_2} & \cdots & \delta_{i_2}^{j_{q-1}}& \delta_{i_2}^{j_q} \\
 \vdots & \vdots & \ddots & \vdots & \vdots\\
 \delta_{i_{q-1}}^{j_1} & \delta_{i_{q-1}}^{j_2} & \cdots & \delta_{i_{q-1}}^{j_{q-1}}& \delta_{i_{q-1}}^{j_q} \\
 \delta_{i_q}^{j_1} & \delta_{i_q}^{j_2} & \cdots & \delta_{i_q}^{j_{q-1}}& \delta_{i_q}^{j_q} \\
 \end{array}
 \right|
 \eqno (20)
$$
     Assume that (19) is true for $r=k$, we only need to show that
     it is also true for $r=k+1$. For $r=k+1$, Using (12) and
     (20), we have
     $$\begin{array}
       {rcl}
       RHS \ \mbox{of}\  (19) &=&
       \displaystyle{\frac{1}{(k+1)!}\sum_{i_1,\cdots,i_{k+1};j_1,\cdots,j_{k+1}}
     \delta_{i_1\cdots i_{k+1} j}^{j_1\cdots j_{k+1} i}}s_{i_1j_1}\cdots
     s_{i_{k+1}j_{k+1}}\\
     &=&
       \displaystyle{\frac{1}{(k+1)!}}\sum
     \left|
                                                      \begin{array}{ccccc}
                                                        \delta_{i_1}^{j_1} & \delta_{i_1}^{j_2} & \cdots & \delta_{i_1}^{j_{k+1}} & \delta_{i_1}^i \\
                                                        \delta_{i_2}^{j_1} & \delta_{i_2}^{j_2} & \cdots & \delta_{i_2}^{j_{k+1}}& \delta_{i_2}^i \\
                                                        \vdots & \vdots & \ddots & \vdots & \vdots\\
                                                        \delta_{i_{k+1}}^{j_1} & \delta_{i_{k+1}}^{j_2} & \cdots & \delta_{i_{k+1}}^{j_{k+1}}& \delta_{i_{k+1}}^i \\
                                                        \delta_j^{j_1} & \delta_j^{j_2} & \cdots & \delta_j^{j_{k+1}}& \delta_j^i \\
                                                      \end{array}
       \right|s_{i_1j_1}\cdots s_{i_{k+1}j_{k+1}}\\
       &=&
       \displaystyle{\frac{1}{(k+1)!}}\sum(\delta_j^i\delta_{i_1\cdots
       i_{k+1}}^{j_1\cdots
       j_{k+1}}-\delta_j^{j_{k+1}}\delta_{i_1\cdots i_ki_{k+1}}^{j_1\cdots j_ki}+\cdots)s_{i_1j_1}\cdots s_{i_{k+1}j_{k+1}}\\
       &=& \sigma_{k+1}\delta_{ij}-
       \displaystyle{\frac{1}{(k+1)!}}\sum\delta_j^{j_{k+1}}\delta_{i_1\cdots i_ki_{k+1}}^{j_1
       \cdots j_ki}s_{i_1j_1}\cdots
       s_{i_{k+1}j_{k+1}}+\cdots\\
       &=& \sigma_{k+1}\delta_{ij}-\sum
       (P_k)_{ii_{k+1}}s_{i_{k+1}j}\\
       &=& (P_{k+1})_{ij}.
     \end{array}$$
   \end{proof}
   \begin{lemma}\label{lemma3.3}
   For each $r$, we have

  $(i)$. $\sum_j(P_r)_{jij}=0$,

   $(ii)$. $\trace(P_rS_F)=(r+1)\sigma_{r+1}$,

   $(iii)$. $\trace(P_r)=(n-r)\sigma_r$.

   \end{lemma}
\begin{proof}

(i). Noting $(j,j_r)$ is skew symmetric in $\delta_{i_1\cdots i_r i}^{j_1\cdots j_r j}$  and $(j,j_r)$ is symmetric in $s_{i_1j_1}\cdots
s_{i_rj_rj}$ (from Lemma 3.1), we have

 $$\sum_{j}(P_r)_{jij}=\displaystyle{\frac{1}{(r-1)!}\sum_{i_1,\cdots,i_r;j_1,\cdots,j_r;j}
     \delta_{i_1\cdots i_r i}^{j_1\cdots j_r j}}s_{i_1j_1}\cdots
     s_{i_rj_rj}=0.
$$

 (ii).  Using (19) and (12), we have
  $$\begin{array}
    {rcl}
    \trace(P_rS_F) &=& \sum_{ij}(P_r)_{ij}s_{ji}\\
    &=&
    \frac{1}{r!}\sum_{i_1,\cdots,i_r;j_1,\cdots,j_r;i,j}\delta_{i_1\cdots
    i_rj}^{j_1\cdots j_ri}s_{i_1j_1}\cdots s_{i_rj_r}s_{ji}\\
    &=& (r+1)\sigma_{r+1}.
  \end{array}
  $$

  (iii). Using (ii) and the definition of $P_r$, we have
  $$\trace(P_r)=\trace(\sigma_rI)-\trace(P_{r-1}S_F)=n\sigma_r-r\sigma_r=(n-r)\sigma_r.$$
  \end{proof}
\noindent {\bf Remark 3.1} When $F=1$, Lemma 3.3 was a well-known result (for example, see Barbosa-Colares [1]).

\begin{lemma}\label{lemma3.5}
  If $\lambda_1=\lambda_2=\cdots=\lambda_n=\mbox{const}\neq0$, then up to translations and
homotheties,
  $X(M)$ is the Wulff shape.
\end{lemma}
\begin{proof}
 Choose a local orthogonal frame $e_1, e_2, \cdots, e_n$ such that $A_F$ is diagonalized:
  $$A_F=diag(\mu_1,\cdots,\mu_n),
\eqno (21)
  $$
  where $\mu_i>0$ for $i=1,\cdots, n$ by the convexity condition.
  Then we have $S_{ij}=\mu_i h_{ij}$. From (10) and (12), we get
  $$\begin{array}{rl}
  &0=M_1^2-M_2=(\frac1n\sum_i\mu_ih_{ii})^2-\frac{2}{n(n-1)}\sum_{i<j}\mu_i\mu_j(h_{ii}h_{jj}-h_{ij}^2)\\
  =&\frac{1}{n^2(n-1)}\{(n-1)(\sum_i\mu_ih_{ii})^2-2n\sum_{i<j}\mu_i\mu_j(h_{ii}h_{jj}-h_{ij}^2)\}\\
  =&\frac{1}{n^2(n-1)}\sum_{i<j}\{(\mu_ih_{ii}-\mu_jh_{jj})^2+2n\mu_i\mu_jh_{ij}^2\},\end{array}$$
  so, $\mu_1h_{11}=\mu_2h_{22}=\cdots=\mu_nh_{nn}$ and $h_{ij}=0$
  when $i\neq j$. Then, from [10] or [3], [14],
up to translations and homotheties,
  $X(M)$ is the Wulff shape.
\end{proof}

\section {Proofs of Theorem 1.1-Theorem 1.5}
{\bf  Proof of Theorem 1.1}
  By use of (5), we have
$$
\langle X, \nu\rangle_i=-\sum_jh_{ij}\langle X, e_j\rangle,\quad
        \langle X, e_j\rangle_i=\delta_{ij}+h_{ij}\langle X, \nu\rangle,
\eqno(22)
$$
so, from (16), Lemma 3.1 and (i), (ii), (iii) of Lemma 3.3, we have the following calculation
$$
\begin{array}
   {rcl}
 &&\DIV\{P_r(\langle X, \nu\rangle\grad_{S^n}F-F\grad|X|^2/2)\}\\
 &=& \sum_{ij}\{(P_r)_{ij}(\langle X,\nu\rangle F_j-F\langle X, e_j\rangle)\}_i\\
 &=& \sum_{ij}(P_r)_{ij}\{-\sum_kh_{ik}(\langle X, e_k\rangle F_j+\langle X, \nu\rangle F_{jk}-F_k\langle X,
 e_j\rangle)
 -F\delta_{ij}-Fh_{ij}\langle X, \nu\rangle\}\\
 &=& -\sum_{ijk}h_{ki}(P_r)_{ij}\langle X, e_k\rangle F_j+\sum_{ijk}h_{ki}(P_r)_{ij}\langle X, e_j\rangle F_k\\
 && -\langle X, \nu\rangle\sum_{ijk}(P_r)_{ij}(F_{jk}+F\delta_{jk})h_{ki}-F\sum_i(P_r)_{ii}\\
  &=& -\sum_{ijk}h_{ki}(P_r)_{ij}\langle X, e_k\rangle F_j+\sum_{ijk}h_{ji}(P_r)_{ik}\langle X, e_k\rangle F_j\\
 && -\langle X, \nu\rangle\sum_{ijk}(P_r)_{ij}A_{jk}h_{ki}-F\sum_i(P_r)_{ii}\\
 &=& -\langle X, \nu\rangle\sum_{ij}(P_r)_{ij}s_{ji}-F\sum_i(P_r)_{ii}\\
 &=& -\langle X, \nu\rangle\trace(P_rS_F)-F\trace(P_r)\\
 &=& -\langle X, \nu\rangle(r+1)\sigma_{r+1}-F(n-r)\sigma_r\\
 &=& -(n-r)C_n^r(FM_r+M_{r+1}\langle X, \nu\rangle).
        \end{array}$$
Integrating the above formula over $M$, we get (2) by use of Stokes
Theorem.

\paragraph{Proof of Theorem 1.2:}
From (2), we have
$$
\int_M(F+M_1\langle X, \nu\rangle)\dif A_X=0, \eqno (23)
$$
$$
\int_M(FM_1+M_2\langle X, \nu\rangle)\dif A_X=0. \eqno (24)
$$

By the assumption $M_1={\rm constant}$, we get from (23) and (24)
$$\int_M\langle X, \nu\rangle(M_1^2-M_2)\dif A_X=0.
\eqno (25)
$$

On the other hand,
$$
M_1^2-M_2=\frac{1}{n^2(n-1)}\sum_{j<i}(\lambda_i-\lambda_j)^2\geq0. \eqno (26)
$$
Thus, if $\langle X, \nu\rangle$ has fixed sign, then $M_1^2-M_2=0$,
so
$$\lambda_1=\lambda_2=\cdots=\lambda_n.$$
 Thus, from
Lemma 3.4, up to translations and homotheties,
  $X(M)$ is the Wulff shape.

\paragraph{Proof of Theorem 1.3:}
We have the fact that if $M$ is compact and $M_r>0$ then
 $$
 \label{mk}
 M_{r-1}\geq M_r^{(r-1)/r},\quad 2\leq r\leq n
 \eqno (27)
 $$
 with equality holds if and only if $\lambda_1=\lambda_2=\cdots=\lambda_n$ on $M$ (cf. [9], [2]).
 Indeed (27) holds if $M_r\equiv\mbox{const}$,  since $M$ is compact, there exists a point $p_0$ on $M$ such that all
 principal curvatures are positive at $p_0$,
 so all anisotropic principal curvatures are positive at $p_0$. Applying (27) inductively, one sees that if $M_r\equiv\mbox{const}$, then
$$
M_r\leq M_1^r, \eqno (28)
$$
here again equality holds if and only if $\lambda_1=\lambda_2=\cdots=\lambda_n$ .

  Integrating $FM_r^{(r-1)/r}\leq FM_{r-1}$ over $M$, using (2) and $M_r={\rm constant}$, we get
$$
M_r^{(r-1)/r}\int_MF\dif A_X\leq \int_MFM_{r-1}\dif A_X=-M_r\int_M\langle X, \nu\rangle\dif A_X. \eqno (29)
$$
  On the other hand, our assumption $M_1={\rm constant}$ (thus $M_1>0$) and (23) implies
  $$
  \int_M\langle X, \nu\rangle\dif A_X=-\frac{1}{M_1}\int_MFdA_X.
  \eqno (30)
  $$

Putting (30) into (29), we get
  $$M_1^r\leq M_r.
\eqno (31)
  $$
  Therefore equality holds in (28) and $\lambda_1=\lambda_2=\cdots=\lambda_n$ on
  $M$. Thus, from Lemma 3.4, up to translations and homotheties,
  $X(M)$ is the Wulff shape.

\paragraph{Proof of Theorem 1.4:}

From (2), we have

$$\int_M(FM_k+M_{k+1}\langle X, \nu\rangle)\dif A_X=0.
\eqno (32)
$$
$$\int_M(FM_r+M_{r+1}\langle X, \nu\rangle)\dif A_X=0.
\eqno (33)
$$
From the assumptions $\frac{M_r}{M_k}={\rm constant}$, $\frac{M_r}{M_k}\times (32)-(33)$ implies
$$
\int_M\langle X, \nu\rangle(M_{r+1}-\frac{M_r}{M_k}M_{k+1})\dif A_X=0. \eqno (34)
$$

From the convexity of $M$, all the principal curvatures of $M$ are positive, so all the anisotropic principal curvature are positive, we have
$M_l>0$, $0\leq l\leq n$ on $M$. From
$$
M_k\cdot M_{k+2}\leq M^2_{k+1},\quad \cdots, \quad M_{r-1}M_{r+1}\leq M^2_r, \eqno (35)
$$
where equality holds in one of (35) if and only if $\lambda_1=\lambda_2=\cdots=\lambda_n$, we can check
$$
M_kM_{r+1}\leq M_{k+1}M_r,
$$
that is
$$
M_{r+1}-\frac{M_r}{M_k}M_{k+1}\leq 0. \eqno (36)
$$
On the other hand, we can choose the position of origin $O$ such that $\langle X, \nu\rangle$ has fixed sign. Thus, from (34) and (36),
$M_kM_{r+1}=M_{k+1}M_r,$ so $\lambda_1=\lambda_2=\cdots=\lambda_n.$
 Thus, from
Lemma 3.4, up to translations and homotheties,
  $X(M)$ is the Wulff shape.

\paragraph{Proof of Theorem 1.5:}
From proof of Theorem 1.3, we know that there exists a point $p_0\in M$ such that the anisotropic principal curvature $\lambda_i(p_0)>0$, $1\leq
i\leq n$. From $\frac{M_k}{M_n}={\rm constant}$, we have $\frac{M_k}{M_n}=\frac{M_k}{M_n}(p_0)>0$. Thus $M_n\not=0$ on $M$, by the continuity of
$\lambda_i$, we have $\lambda_i>0$, $1\leq i\leq n$, on $M$. Therefore all principal curvatures of $M$ are positive on $M$ and $M$ is convex.
Theorem 1.5 follows from Theorem 1.4.


\vskip 1cm
\begin{flushleft}
Department of Mathematical Sciences\\
Tsinghua University\\
Beijing 100084\\
P. R. China\\
\texttt{yjhe@math.tsinghua.edu.cn} \\
\texttt{hli@math.tsinghua.edu.cn}
\end{flushleft}
\end{document}